\documentclass[journal]{IEEEtran}
\usepackage{multirow}
\usepackage{comment}
\usepackage{soul}
\usepackage{color}
\usepackage[usenames,dvipsnames]{xcolor}

\usepackage{cite}
\usepackage{url}

\ifCLASSINFOpdf
\else
\fi
\usepackage{tikz}
\usetikzlibrary{arrows, positioning, shapes.geometric, fit, calc,backgrounds,shapes.arrows}
\usetikzlibrary{arrows.meta}

\usepackage{amsmath}
\usepackage{amsfonts}
\usepackage[short]{optidef}
\usepackage[ruled,vlined]{algorithm2e}
\usepackage{ bbold }
\usepackage{amssymb}

\ifCLASSOPTIONcompsoc
  \usepackage[caption=false,font=normalsize,labelfont=sf,textfont=sf]{subfig}
\else
  \usepackage[caption=false,font=footnotesize]{subfig}
\fi

\usepackage{flushend}

\usepackage{bm}
\usepackage[breaklinks]{hyperref}
\usepackage[capitalise]{cleveref}
\Crefname{equation}{}{}

\begin{document}



\title{Network-Aware Flexibility Requests for Distribution-Level Flexibility Markets}

\author{Eléa~Prat, Irena~Dukovska, Rahul~Nellikkath, Malte~Thoma, Lars~Herre, Spyros~Chatzivasileiadis
  \thanks{Eléa Prat, Lars Herre, Rahul Nellikkath, and Spyros Chatzivasileiadis are with the Technical University of Denmark, 2800 Kgs. Lyngby, Denmark. e-mail: \{emapr, lfihe, rnelli, spchatz\}@dtu.dk. Irena Dukovska is with Technical University of Eindhoven. e-mail: i.dukovska@tue.nl. Malte Thoma is with bnNETZE, email: malte.thoma@bnnetze.de. This work is supported by the H2020 European Project FLEXGRID, Grant Agreement No. 863876.}
}

\maketitle

\IEEEpubidadjcol

\begin{abstract}
This paper proposes a method to design network-aware flexibility requests for local flexibility markets. These markets are becoming increasingly important for distribution system operators (DSOs) to ensure grid safety while minimizing costs and public opposition to new network investments. Despite extended recent literature on local flexibility markets, little attention has been paid to quantifying the flexibility required at each location, considering physical network constraints (e.g. line and voltage limits). The method introduced uses a chance-constrained optimization model and a \emph{LinDistFlow} approximation to consider both physical network constraints and uncertainty caused by renewable production or demand fluctuations. Unlike other methods, it avoids sharing sensitive grid data with the market operator. We compare our approach against a stochastic market-clearing mechanism which serves as a benchmark, and we derive analytical conditions for the performance of our method to determine flexibility requests. We show on two case studies that our method outperforms the stochastic market-clearing benchmark in terms of computation time while achieving comparable social welfare and costs for the DSOs. One of the case studies is conducted on an actual German distribution grid, showing that the proposed method can scale well to real-sized networks.

\end{abstract}


\IEEEpeerreviewmaketitle

\section{Introduction} \label{sec:intro}

Local flexibility markets can help defer or even avoid significant distribution network investments by tapping into the available flexibility of local resources (e.g. home batteries, heating and cooling systems). Considering the often strong public opposition against new network infrastructure, local flexibility markets are expected to play a significant role in the efforts of the distribution system operators (DSOs) to maintain a safe and reliable grid operation. In that respect, the DSO is expected to be one of the key flexibility market players, acting as a flexibility buyer, in order to ensure grid security in case of contingencies and generation uncertainties, similar to the role of transmission system operators in the transmission-level ancillary services market.
Here, we define flexibility as a controlled real-time power variation from a pre-defined setpoint. Flexibility can be traded in energy markets, which settle activation, or in capacity markets, with a guarantee on availability. In this paper, we consider flexibility trading in capacity markets.

The concept of \textit{local flexibility markets} at the distribution level has emerged recently \cite{Jin2020}, and is investigated in several EU projects\footnote{See \url{http://www.interrface.eu/}, \url{http://smartnet-project.eu/}, \mbox{\url{https://flexgrid-project.eu/}}, and \url{http://www.eu-ecogrid.net/}.}, such as INTERRFACE, SmartNet, and FLEXGRID, as well as national projects such as EcoGrid 2.0 in Denmark. A large body of work has recently focused on local flexibility markets, including bidding strategies for flexibility offers \cite{Olivella-Rosell2018}, and efforts to account for network constraints \cite{Prat2021,Fonteijn2020}, uncertainty \cite{Garcia-Torres2021}, and co-optimization of energy and reserves~\cite{Bouffard2005}.

Considering the continuously increasing electricity demand, bidirectional flows, and reduced investments in network infrastructure, the distribution network can no longer be considered a copperplate \cite{FLEXGRID}. As a result, local flexibility markets are expected to implicitly or explicitly consider the network constraints in the market clearing (e.g. \cite{Prat2021}), so that the procured flexibility helps resolve and not further aggravate existing network problems. Here, we qualify as ``network-aware'' any approach in which the distribution network is not considered to be a copperplate.

For network-aware local flexibility markets, the DSO has two options: information sharing or privacy preserving. The first option is to share all network data with the flexibility market operator (FMO), as proposed in \cite{Ullah2021,Olivella-Rosell2018,Khorasany2021,Fonteijn2020}. 
In \cite{Papavasiliou2018}, the authors move even a step further and present an approach toward understanding distribution locational marginal prices by decomposing the distribution locational marginal price (DLMP) for \emph{energy} into terms relating to power at the root node, to real power losses, to reactive power losses, to voltage constraints, and to line limits.

This information-sharing solution has received criticism for (i) the potentially intractable and impractical size of the problem for real-life operation, especially when considering a stochastic framework, and (ii) due to the existing legal framework on data exchange between DSO and FMO. Another option would be for the DSO to take over the role of the FMO. 
However, the role of the FMO is currently not addressed in the EU framework \cite{roadmap_flex}, namely the Regulation 2019/943 \cite{ERegulation} and the Directive  2019/944 \cite{EDirective} for the internal market for electricity. Therefore, the implementation of such an approach is still unclear and may differ per country, while it may not be practically implementable in some countries. Besides, having DSO and FMO as separate entities is the most common set-up in the literature, as shown in the review \cite{Jin2020}, and it is chosen by many pilot projects on flexibility markets, which is discussed in \cite{schittekatte2020flexibility}. For these reasons, we make the assumption that DSO and FMO are different actors.

The second, the privacy-preserving option, entails the creation of a network-aware \emph{FlexRequest} by the DSO, which is then used in a deterministic market clearing by the FMO, similar to the transmission-level ancillary services markets currently operating in Europe. This privacy-preserving option is the focus of this paper.
It is compatible with real-world operations and ensures high transparency for electricity traders since the FMO does not need to consider stochasticity and network constraints in the market clearing. Thus, it is compatible with the requirements of the EU Directive \cite{EDirective} for transparent and standard market products for flexibility services. However, the question is: How to determine such a network constrained \emph{FlexRequest} under uncertainty?

Many existing works on flexibility markets disregard the problem of how DSOs shall estimate how much flexibility they need. In \cite{Prat2021,Fonteijn2021,Ullah2021,Olivella-Rosell2018,Garcia-Torres2021,Zhao2021}, flexibility requests from the DSO are treated as purely external input parameters. Their focus is instead on the optimization problem of the market operator or the balance responsible party.
The framework of \cite{Khorasany2021} provides a heuristic approach for determining flexibility requests which, however, assumes a priori knowledge of each prosumer's power exchange, and does not consider network constraints. Within the flexibility mechanism of \cite{Fonteijn2020}, an approach for obtaining the time and value of needed flexibility (similar to \emph{FlexRequest}) is proposed for transformer overloading, based on the cost of transformer lifetime reduction. However, for line congestions and voltage support across the network, there exists no generic real-world approach to determine flexibility requests. In general, capacity markets for flexibility are considerably less investigated than energy markets for flexibility. Among the few works on capacity flexibility markets is \cite{Prat2021}, which, however, does not address the problem of \emph{FlexRequest} creation.

In the studies on energy flexibility markets, deterministic formulations are implemented. Mechanism design is used to design a fair market for energy flexibility at the distribution level in \cite{tsaousoglou2021}. A deterministic convex Second Order Conic Program (SOCP) is used to model the AC power flow and determine the required flexibility. The coordination of energy flexibility markets in the distribution and transmission system is studied in \cite{marques2023}, where the linearized \emph{LinDistFlow} model is used for modeling the distribution network. The uncertainty of flexibility a DSO provides to the point of common coupling stemming from PV and battery systems is modeled in \cite{rayati2022}. However, the flexibility request is considered as an external input. A model for coordination of a flexibility market for energy and flexibility products with a single FMO and several DSOs is modeled in \cite{paredes2023}. Uncertainty is incorporated in the modeling of the flexibility offers to the market, not the request for flexibility.

To the best of our knowledge, this is the first paper that proposes a framework to design network-aware \emph{FlexRequests} at the distribution level. The contributions of this paper are:
\begin{itemize}
  \item We propose a tool for the DSO to determine network-aware \emph{FlexRequests} that will be submitted in a distribution-level capacity flexibility market, run by an independent FMO. This tool also considers uncertainty and preserves the privacy of both the DSO and flexibility providers.
  \item We study the trade-off between the privacy-preserving property ensured by our method and ideal flexibility procurement obtained through stochastic market clearing, by deriving the analytical conditions under which our proposed method matches the performance of stochastic market clearing.
  \item We demonstrate the use of this new tool coupled with a deterministic market clearing and we examine the influence of market liquidity and the definition of bidding zones by the DSO. We compare this setup to the stochastic market-clearing benchmark. We show how our approach outperforms the stochastic market clearing in terms of computation time while achieving similar social welfare and costs for the DSO. Using a real-world German distribution network, we show that our method scales better than the benchmark.
\end{itemize}

The remainder of this paper is organized as follows. Section~\ref{sec:framework} provides details on the market framework used in the rest of the paper. The \emph{FlexRequest} creation optimization problem is presented in Section~\ref{sec:Mcreation}. Both the deterministic and stochastic market-clearing mechanisms are outlined in Section~\ref{sec:marketclearing}. We evaluate the suboptimality of our approach compared to the stochastic market clearing in Section~\ref{sec:subopt}. The case studies and the corresponding results are given in Section~\ref{sec:Case}. We draw conclusions in Section~\ref{sec:Conc}.

\section{General Framework} \label{sec:framework}

\begin{figure}[th]
    \centering
    \includegraphics[width=\linewidth]{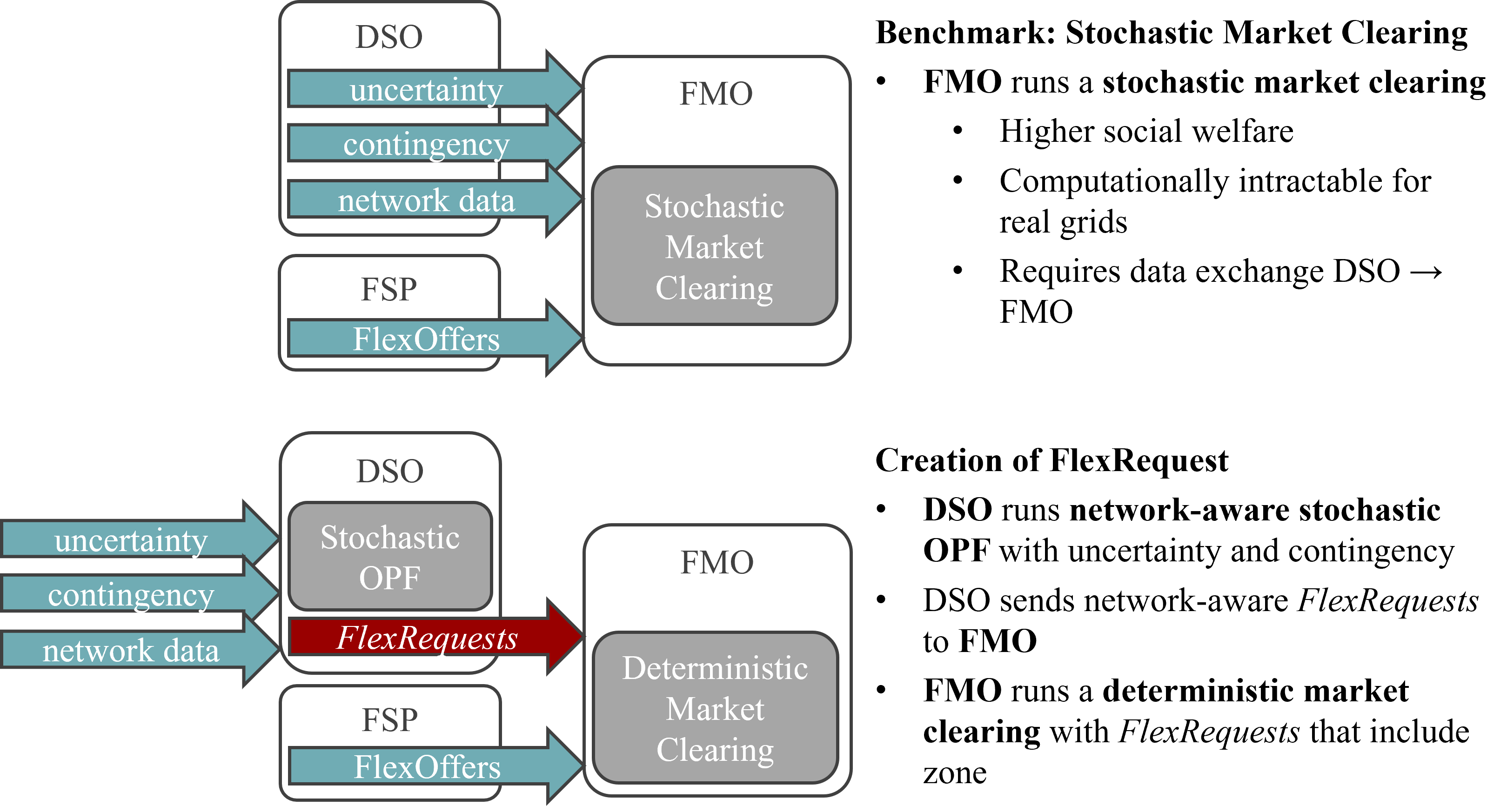}
    \caption{Conceptual representation of the proposed market design. Top: Benchmark option with a stochastic market clearing where the FMO requires network data. Bottom: Proposed privacy-preserving \emph{FlexRequest} market design where the FMO only requires knowledge of the zone (DSO area).}
    \label{fig:Overview}
    \vspace{-3mm}
\end{figure}

As mentioned in the introduction, there are two options for flexibility markets, depending on the desire or possibility for the DSO to share its network data with the FMO. These are represented in Figure \ref{fig:Overview}. The focus is on flexibility markets in which the DSO is the single buyer and the procured flexibility (reserves) are used to solve expected problems in the network. Our approach explicitly considers the underlying network and its associated constraints (line congestion, voltage limits), as well as the uncertainty of fluctuating renewable production or loads (e.g. electric vehicles).

The first setup is a stochastic market clearing or security-constrained unit commitment as in, e.g., [6]. In this approach, the FMO finds the optimal solution to the social welfare maximization problem and explicitly determines the quantity and location of reserves. The resulting stochastic problem is computationally intensive, which may make it impractical and even intractable for real-sized distribution networks. For this approach, information on the distribution network, as well as on expected uncertainty and/or contingency must be communicated to the FMO by the DSO. This data is usually private to the DSO.
Furthermore, stochastic programming also requires a complete market design overhaul \cite{Papakonstantinou2016} which is why the few real-world examples of stochastic electricity markets are limited to exogenous sizing of probabilistic security margins (reserves) within otherwise deterministic market-clearing routines \cite{Dvorkin2020}.

The second setup, presented in the bottom part of Figure \ref{fig:Overview} is a privacy-preserving deterministic reserve flexibility market setup. Compared to the other setup, this approach does not require that private data be communicated to other parties by the DSO. The DSO determines network-aware \emph{FlexRequests} considering uncertainty, by solving a stochastic optimal power flow (OPF) problem. The drawback of this method is that it may not achieve the highest social welfare, though it may come close to it. The FMO in this case solves a deterministic market model that is transparent for all participants and computationally less intensive. As such, it can be easily integrated into the existing market frameworks in Europe. In \cref{sec:subopt}, we analytically quantify the tradeoffs in terms of social welfare in more detail.

In the proposed market framework, the DSO takes into account the uncertainty in the power injections when creating \emph{FlexRequests}.
We consider that the DSO's objective is not to cover imbalances, which are the responsibility of the balance responsible parties. The aim of the \emph{FlexRequests} is rather to avoid line congestion and violation of voltage limits. In case of imbalances resulting from uncertainty realization, we assume that the difference is covered by exchange with the upper-level grid, through the slack bus.
The DSO is expected to know its network sufficiently well to define zones, i.e., DSO areas. Something similar happens in transmission-level markets, where the zones are defined primarily with the help of the transmission system operators. Following that, the DSO submits the corresponding \emph{FlexRequest} volume and price to the FMO, along with the zones defined. Defining zones on a short-term basis and re-evaluating them frequently allows for more flexibility. As the markets are local, zones have the capability to change more often than at the transmission level, depending on the local assets and their consumption and production patterns. At the same time, the flexibility service providers (FSP) submit their \emph{FlexOffers} to the FMO. The FSPs are distributed resources located in the distribution network, but no further assumptions are made about their properties, which can be diverse. The decision-making process to determine \emph{FlexOffers}\footnote{The decision-making process to determine \emph{FlexOffers} is investigated in, e.g., \mbox{\url{https://flexgrid-project.eu/}}}~\cite{Shen2020,Olivella-Rosell2020} is out of the scope of this paper.

The market-clearing mechanism that is used for clearing the \emph{FlexRequests} (\cref{sec:Mdeter}) is kept very general so that it can be adapted to different frameworks. In particular, we do not impose assumptions on the optimization horizon, time resolution, pricing, or payments.

This paper focuses on flexibility markets for active power only, as these are considered to have the highest value at the moment. Congestion problems are becoming a pressing issue for DSOs, and flexibility markets are seen as an alternative or supporting solution to network reinforcements for solving these problems. Therefore, there is an extensive body of literature addressing the design of flexibility markets focusing predominantly on markets for active power. Nevertheless, in our work, we also include reactive power, assuming that it is linked to the active power by a constant factor.
This allows for indirect consideration of reactive power in the active power flexibility market, without the need to resort to complex markets for active and reactive power, for the moment. The extension to a market that explicitly considers reactive power flexibility is part of future work and is out of the scope of this study. Finally, we focus on radial distribution networks, which are the most common topology used in distribution systems.

\section{Creation of FlexRequest by the DSO} 
\label{sec:Mcreation}

As mentioned in the introduction, the need for the DSO to procure flexibility is motivated by the need to avoid investments. Normally, the DSO shall guarantee the safe/secure operation of the distribution grid under a wide range of possible conditions. Considering the increasing uncertainty in supply and demand at the distribution level, this need for investments becomes prohibitively high. Instead, the DSO can procure cheaper flexibility reserves from distributed energy resources (DER) and loads located in the distribution network through a market mechanism. The reserves are procured beforehand and will be (partially) activated in real time depending on the realization of the uncertainty. The head time between reserves procurement and real-time activation can vary from several hours ahead to a day ahead, depending on the design of the market. The DSO does not own or control the resources that can provide flexibility. Therefore, it has to estimate the quantity and direction of the required flexibility in the form of \emph{FlexRequests} that will be sent to the FMO for market clearing. 

In this section, we present how the DSO models the uncertainty stemming from renewable generation using chance-constrained optimization. We explain the motivation for the chosen framework and introduce the equations of the model. We further present a deterministic reformulation, which is required for tractability. Finally, we address the questions of defining zones in the flexibility market and adequate pricing of the \emph{FlexRequests}.

\subsection{Motivations for Modelling Choices}
To incorporate uncertainty in the modeling, in this paper, we chose to apply chance-constrained optimization\cite{Geng2019}; we motivate our choices below. 

The chance-constrained modeling framework is appropriate for modeling network-constrained problems in which the constraints have to be satisfied within some probability level. 
Since line limits can, in practice, occasionally be exceeded, chance-constrained optimization is a well-suited modeling approach for creating \emph{FlexRequests}. Moreover, it allows the DSOs to decide the acceptable violation probability per asset or per constraint \cite{Roald2018}. Allowing for a small probability of violation provides a risk-aware and cost-effective mechanism as compared to conservative, and thus expensive solutions obtained with robust optimization \cite{Ning2019}. Through the implementation of linear policies, the relation between the uncertainty realization and the system reaction can be modeled. In this way, a response to scenarios that are not represented in the scenario set is ensured, which is not guaranteed when applying stochastic optimization with scenarios. 

The approach for the creation of \emph{FlexRequests} is somewhat similar to the dispatch of controllable generators under uncertainty, but with the added complexity that the DSO does not have control over the resources providing the reserves.
For the dispatch problem, a distributionally robust chance-constrained OPF for distribution systems is presented in \cite{Mieth2018} and \cite{Dall'Anese2017}, using the \emph{LinDistFlow} approximation \cite{baran1989} to model the network. The only uncertainty considered in \cite{Mieth2018, Dall'Anese2017}, is in the voltage limits, whereas line limits are additionally considered in \cite{Mieth2020}.
Existing studies that consider chance-constrained OPF in the distribution network, assume that the DSO can fully control the resources required to compensate for the effects of realizations of uncertainty via direct dispatch and that they provide certain responses.

A case in which the network operators do not own the resources that provide flexibility and their response is uncertain, such as the case covered in this study, is considered in \cite{Zhang2017}. There, uncertain aggregated loads can provide market balancing reserves, in addition to controllable generators.
However, \cite{Zhang2017} focuses on transmission networks and develops its approach with a DC linearized power flow, which is not appropriate for distribution networks. 

To the best of our knowledge, the design of \emph{FlexRequests} and the procurement of flexibility reserves by the DSO through market mechanisms to solve uncertain line congestions and voltage problems in the distribution network is still an open topic that needs to be addressed.

\subsection{Chance-Constrained Model} \label{CC_model}

The objective of the DSO is to ensure that the operational network constraints are maintained with a given probability. The objective function~\eqref{OF_flexreq} is formulated as a minimization of the total required upward and downward quantity of flexibility requests which will be sufficient to respond to uncertainty realizations in real time.

The optimization problem for the DSO is formulated as
\begin{subequations}
\begingroup
\allowdisplaybreaks
\begin{align}\label{OF_flexreq}
     \min_\mathbf{x} &\sum_{t \in \mathcal{T}} \sum_{n \in \mathcal{N}} (P_{n,t}^\text{R+} + P_{n,t}^\text{R-}) \\
     \text{s.t.} 
     \label{active_p} 
     \nonumber & \sum_{ij \in \mathcal{L}, j=n} P_{ij,t} + P_{n,t}^\text{inj}+P_{n,t}^\text{R}      \\
     &  \quad + P^\text{U}_{n, t} = \sum_{ij \in \mathcal{L}, i=n} P_{ij,t}, && \forall n \in \mathcal{N}, \ \forall t \in \mathcal{T}, \\
     \nonumber & \sum_{ij \in \mathcal{L}, j=n} Q_{ij,t} + Q_{n,t}^\text{inj}+Q_{n,t}^\text{R}    \\  
     & \quad + Q^\text{U}_{n,t} = \sum_{ij \in \mathcal{L}, i=n} Q_{ij,t}, && \forall n \in \mathcal{N}, \ \forall t \in \mathcal{T}, \label{reactive_p} \\
    \nonumber & u_{j,t} =  u_{i,t}  \\
    & \quad - 2(R_{ij}P_{ij,t}+X_{ij}Q_{ij,t}), && \forall ij \in \mathcal{L}, \ \forall t \in \mathcal{T},  \label{sq_voltage} \\
    & \tilde{P}_{n,t}^\text{R}(\xi) =  P_{n,t}^\text{R} + \alpha_{n,t} \xi_{\text{tot},t}^\text{P}, && \forall n \in \mathcal{N}, \ \forall t \in \mathcal{T}, \label{part_factor} \\
    &\nonumber - \sum_{ij \in \mathcal{L}, i=n} P_{ij,t}  +  P_{n,t}^\text{inj}  + \\ 
    & + \tilde{P}_{n, t}^R - \sum_{k \in \mathcal{N}, k \neq \text{ref}} \tilde{P}_{k,t}^\text{R} =0, &&  n=\text{ref}, \ \forall t \in \mathcal{T}, \label{pref} \\
    & \nonumber- \sum_{ij \in \mathcal{L}, i=n} Q_{ij,t} + Q_{n,t}^\text{inj}  + \\ 
    & + \tilde{Q}_{n, t}^R - \sum_{k \in \mathcal{N}, k \neq \text{ref}} \tilde{Q}_{k,t}^\text{R} = 0, &&  n=\text{ref}, \ \forall t \in \mathcal{T}, \label{qref}\\
    & \sum_{n \in \mathcal{N}, n \neq \text{ref}} \alpha_{n,t} = 0, && \forall t \in \mathcal{T}, \label{sum_bal2} \\
    & P_{n,t}^\text{R-}, P_{n,t}^\text{R+} \geq 0, && \forall n \in \mathcal{N}, \ \forall t \in \mathcal{T}, \label{req_def} \\
    \nonumber & \mathbb{P}(\tilde{P}_{ij,t}^2(\xi) + \tilde{Q}_{ij,t}^2(\xi) \leq \overline{S}_{ij}^2) \\
    & \quad  \geq 1 - \epsilon_S, && \forall ij \in \mathcal{L}, \ \forall t \in \mathcal{T} \label{cc_apparent}\\
    & \mathbb{P}(\underline{V}_n^2 \leq \tilde{u}_{n,t}(\xi)) \geq 1 - \epsilon_V,  && \forall n \in \mathcal{N}, \ \forall t \in \mathcal{T}, \label{cc_v_down} \\
    & \mathbb{P}(\tilde{u}_{n,t}(\xi) \leq \overline{V}_n^2) \geq 1 - \epsilon_V, &&  \forall  n \in \mathcal{N}, \ \forall t \in \mathcal{T}, \label{cc_v_up} \\
    \nonumber & \mathbb{P}(- P_{n,t}^\text{R-} \leq \tilde{P}_{n,t}^\text{R}(\xi)) \\
    & \quad \geq 1 - \epsilon_R,  && \forall n \in \mathcal{N}, \ \forall t \in \mathcal{T}, \label{cc_req_down} \\
    & \mathbb{P}(\tilde{P}_{n,t}^\text{R}(\xi) \leq P_{n,t}^\text{R+}) \geq 1 - \epsilon_R,  && \forall n \in \mathcal{N}, \ \forall t \in \mathcal{T}, \label{cc_req_up} 
\end{align}
\endgroup
\end{subequations}
where $\mathbf{x}=\{P^\text{R+},P^\text{R-},P^\text{R}, u, P, Q\}$

To represent the network, with sufficient accuracy, despite neglecting power losses, the \mbox{\emph{LinDistFlow}} approximation \cite{baran1989} is used in \eqref{active_p}-\eqref{sq_voltage}. This approximation is valid for radial networks. The active and the reactive power balance per node are expressed in \eqref{active_p} and \eqref{reactive_p}. Here, $P_{ij,t}$ and $Q_{ij,t}$ are the active and reactive power flow in the line between bus $i$ and $j$ at time $t$, in case of perfect forecast. Each line $ij$ has a resistance $R_{ij}$ and reactance $X_{ij}$. The sets $\mathcal{N}$, $\mathcal{L}$, and $\mathcal{T}$ gather, respectively, the buses of the system, the lines (expressed with their origin and destination bus), and the time periods considered. The power injection at node $n$ is separated between a deterministic component $P_{n,t}^\text{inj}$, a forecast of the uncertain part $P_{n,t}^\text{U}$ and the flexibility activation $P_{n,t}^\text{R}$, which is a decision variable.
The associated reactive components are $Q_{n,t}^\text{inj}$, $Q_{n,t}^\text{U}$ and $Q_{n,t}^\text{R}$. It is assumed that the reactive power injections are proportional to the active power injections according to the power factor $\cos  \phi $ and through the parameter $K = \sqrt{\frac{1-\cos \phi^2}{\cos \phi^2}}$. The voltage drop is given in \eqref{sq_voltage}. The variable $u_{n,t}$ is an auxiliary variable, which is equal to the squared voltage magnitude. 

The uncertain power injection is defined as: \mbox{$\tilde{P}^\text{U}_{n, t}(\xi) = P^\text{U}_{n, t} - \xi_{n,t}, \ \forall n \in \mathcal{N}, \forall t \in \mathcal{T} $}, where $P^\text{U}_{n, t}$ is the forecasted value and $\xi_{n,t}$ is the deviation from the forecast at bus $n$, for time $t$. The activated flexibility due to the \emph{FlexRequests} $\tilde{P}^R_{n,t}(\xi)$ is expressed as an affine function of the total forecast error $ \xi_{\text{tot},t} = \sum_{n \in \mathcal{N}} \xi_{n,t}, \forall t \in \mathcal{T} $ through the factors $\alpha_n, \ \forall n \in \mathcal{N}$, as in \eqref{part_factor}.

The realizations of uncertainty and the activation of \emph{FlexRequests} are causes for imbalance in the nodal active and reactive balance equations. Since power balancing is the responsibility of the retail company and not the DSO, it is considered that the DSO is not in charge of resolving imbalances. The DSO is responsible for relieving congestions, keeping nodal voltages within limits, and ensuring a safe grid operation. Therefore, the energy necessary to cover them is assumed to be generated in other parts of the network and received through the slack bus. \emph{FlexRequests} are used to ensure that the distribution network will be able to handle the associated changes in terms of line flows and voltages. Since the activation of the \emph{FlexRequests} does not modify the balance, the sum of $\alpha_{n,t}$ of the non-reference buses should be equal to 0, as in \eqref{sum_bal2}. The nodal balance at the slack bus should also be updated accordingly as in \eqref{pref}-\eqref{qref}. In the case of a radial network, the balance is ensured for any realization of the uncertainty when including these constraints. The main variables of the problem, that are submitted to the flexibility market, are $P_{n,t}^\text{R-}$ and $P_{n,t}^\text{R+}$, the downward and upward \emph{FlexRequest} at each bus defined in \eqref{req_def}.

In addition to the flexibility requested, all variables of the problem depend on the uncertainty realization. Following this, five chance constraints are formulated \eqref{cc_apparent}-\eqref{cc_req_up}. Equation \eqref{cc_apparent} is related to the rated apparent power of the lines $\overline{S}_{ij}$. Constraints \eqref{cc_v_down}-\eqref{cc_v_up} require bus voltages to be within limits, between the minimum $\underline{V}_n^2$ and the maximum $\overline{V}_n^2$. Finally in \eqref{cc_req_down}-\eqref{cc_req_up}, the predicted flexibility activation is bounded by the flexibility requested.
The parameters $\epsilon \in (0,1)$ give the violation probability for the corresponding constraint.

The main advantage of this approach for the creation of the \emph{FlexRequests} is that it avoids making assumptions about the flexibility providers. In particular, there is no need to know their position in the network or their flexibility capabilities. In general, the DSO does not have access to such information. For this reason, we assume in \eqref{cc_req_down}-\eqref{cc_req_up} that all the flexibility needs of the DSO will be covered by the flexibility providers. However, if the DSO did have some info on the availability of flexibility, it could be included as supplementary constraints bounding the flexibility activation at the corresponding nodes. The resulting \emph{FlexRequests} would then be more accurate.
Indeed, with the current formulation, there could be multiple equally good solutions in terms of location of the \emph{FlexRequests}, and one of them will randomly be returned by the solver. On the other hand, instead of using a random combination for the locations of the \emph{FlexRequests}, modeling flexibility sources could help select between those solutions by reducing the feasibility space and the number of optimal solutions.

\subsection{Deterministic Reformulation} \label{Det_reformulation}

It is assumed that the forecast error follows a Gaussian probability distribution function with zero mean $\mu_n=0$ and covariance $\Sigma_n$, $\xi_n\sim \mathcal{N}(0, \Sigma_n)$. Given a linear relation between the error and the variables, the chance constraints \eqref{cc_apparent}-\eqref{cc_req_up} can be reformulated analytically to deterministic constraints. Chance constraints of the format $\mathbb{P}(x_i + b_i \xi \leq \overline{x}_i)$ can be reformulated as $x_i + \Phi^{-1}(1-\epsilon_x) \sqrt{b_{i} ^\intercal \Sigma b_{i}} \leq \overline{x}_i$, where $b$ is the matrix that linearly relates the uncertainty source with the variables. Hence, an \emph{uncertainty margin} can be introduced, defined as: $ \Omega_i = \Phi^{-1}(1-\epsilon_x) \sqrt{b_{i} ^\intercal \Sigma b_{i}}$.

Therefore, the linear chance constraints \eqref{cc_v_down}-\eqref{cc_req_up} can be reformulated using uncertainty margins as in \eqref{ref_u}-\eqref{ref_R}:
\begin{subequations} \label{lin_cc}
\allowdisplaybreaks
\begin{align}
    & \underline{V}_n^2 + \Omega^{u}_{n,t} \leq u_{n,t} \leq  \overline{V}_n^2 - \Omega^{u}_{n,t},   & \forall n \in \mathcal{N}, \ \forall t \in \mathcal{T}, \label{ref_u}\\
    \nonumber &-P_{n,t}^\text{R-} + \Omega^{F}_{n,t}  \leq P_{n,t}^\text{R} \\
    & \quad \leq P_{n,t}^\text{R+} -  \Omega^{F}_{n,t} ,  & \forall n \in \mathcal{N}, \ \forall t \in \mathcal{T}. \label{ref_R}
\end{align}
\end{subequations}

A matrix $A$, of size $|\mathcal{L}|\times|\mathcal{N}|$ is introduced, which captures the linear relation between the uncertainty injections and the power flows in the lines \cite{Hassan2018}. Its elements are defined in \eqref{matrix_a}:
\begin{align} \label{matrix_a}
    a_{ij, n}=\begin{cases}
    1,  \text{if $ij$ is on the path from slack bus to $n$}.\\
    0,  \text{otherwise}.
  \end{cases} 
\end{align}

The power flow in the lines \eqref{active_p_lin}-\eqref{reactive_p_lin} and the nodal voltage \eqref{voltage_lin} can be linearly related to the forecast error \cite{Mieth2018}:
\begin{subequations}
\allowdisplaybreaks
\begin{align} 
    \nonumber & \tilde{P}_{ij,t}(\xi) = P_{ij,t} + \sum_{n \in \mathcal{N}}a_{ij, n} (\xi_{n,t}^P-\alpha_{n,t} \xi_{\text{tot},t}^\text{P}), \\
    & \qquad \qquad \qquad \qquad \qquad \qquad \qquad \quad \forall ij \in \mathcal{L}, \ \forall t \in \mathcal{T}, \label{active_p_lin}\\
    \nonumber  & \tilde{Q}_{ij,t}(\xi) = Q_{ij,t} +  \sum_{n \in \mathcal{N}} a_{ij, n} (\xi_{n,t}^\text{Q}-\alpha_{n, t} \xi_{\text{tot},t}^\text{Q}), \\
    & \qquad \qquad \qquad \qquad \qquad \qquad \qquad \quad \forall ij \in \mathcal{L}, \ \forall t \in \mathcal{T} ,\label{reactive_p_lin} \\
    & \tilde{u}_{n,t}(\xi) = u_{n,t} - 2\sum_{ij \in \mathcal{L}}a_{ij, n} \sum_{m \in \mathcal{N}}[R_{ij}a_{ij, m}(\xi_{m, t}^\text{P}-\alpha_{m, t} \xi_{\text{tot},t}^\text{P}) \nonumber \\
    &  + X_{ij}a_{ij, m} (\xi_{m,t}^\text{Q}-\alpha_{m, t} \xi_{\text{tot},t}^\text{Q}) ], \quad \quad \forall n \in \mathcal{N}, \ \forall t \in \mathcal{T}, \label{voltage_lin}
\end{align}
\end{subequations}
where $\Gamma$ is an incidence matrix of size $|\mathcal{N}|\times|\mathcal{U}|$, that denotes the connection of each source of uncertainty to the corresponding bus, where $\mathcal{U}$ is the set that gathers all sources of uncertainty. 

The quadratic chance constraint \eqref{cc_apparent} can be reformulated following the approach of using absolute value chance constraints as in \cite{Lubin2019, Halilbasic2019}. Two auxiliary variables $ k^\text{P} \ \text{and} \  k^\text{Q} $ are introduced and their connection to the rated apparent power of the lines is introduced in \eqref{ref_quad}. The resulting constraints with corresponding uncertainty margins, including the auxiliary variables are defined in \eqref{abs_p}-\eqref{ref_aux_q} respectively: 
\begin{subequations} \label{quad_cc}
\allowdisplaybreaks
\begin{align}
    & (k_{ij,t}^\text{P})^2 +  (k_{ij,t}^\text{Q})^2 \leq \overline{S}_{ij}^2,  && \forall ij \in \mathcal{L}, \ \forall t \in \mathcal{T}, \label{ref_quad}\\
    \nonumber & - k_{ij,t}^\text{P} + \Omega^P_{ij,t} \leq P_{ij,t} \\
    & \quad \leq k_{ij,t}^\text{P} - \Omega^P_{ij,t},  && \forall ij \in \mathcal{L}, \ \forall t \in \mathcal{T}, \label{abs_p}\\
    \nonumber & - k_{ij,t}^\text{Q} + \Omega^Q_{ij,t} \leq Q_{ij,t}  \\
    & \quad \leq k_{ij,t}^\text{Q} - \Omega^Q_{ij,t}, && \forall ij \in \mathcal{L}, \ \forall t \in \mathcal{T}, \label{abs_q} \\
    & k_{ij}^\text{P} \geq \Omega^{k^\text{P}}_{ij,t},  && \forall ij \in \mathcal{L}, \ \forall t \in \mathcal{T}, \label{ref_aux_p}\\
    & k_{ij}^\text{Q} \geq \Omega^{k^\text{Q}}_{ij,t},  && \forall ij \in \mathcal{L}, \ \forall t \in \mathcal{T}. \label{ref_aux_q}
\end{align}
\end{subequations}

After reformulating all the chance constraints, we obtain the following second-order cone (SOC) OPF:
\begin{subequations}
\begin{align}\label{OF_flexreq_ref}
     \min_\mathbf{x} &\sum_{t \in \mathcal{T}} \sum_{n \in \mathcal{N}} P_{n,t}^\text{R+} + P_{n,t}^\text{R-} \\
    \text{s.t.} &\  \text{Eq.} \eqref{active_p}-\eqref{req_def},  \eqref{lin_cc}, \eqref{quad_cc}
\end{align}
\end{subequations}
with $\mathbf{x}=\{P^\text{R+},P^\text{R-},P^\text{R}, u, P, Q, b, \Omega, k^\text{P}, k^\text{Q}\}$

\subsection{Definition of Zones} \label{Zone definition}
With this model, the \emph{FlexRequests} formulated are associated with a given bus. However, in the absence of congestions, the flexibility service could be arbitrarily provided from various buses. The DSO can use the information on potentially congested lines to define zones for the market clearing, as this will increase the liquidity in the market. That is, instead of requesting a flexibility service from a specific bus, the flexibility can be offered by any bus within the same zone. The definition of such zones is out of the scope of this paper, but one idea could be to adapt the approach of \cite{viafora2020}. The determination of dynamic reserve zones is also addressed in~\cite{Wang2015}. In the case studies of Section \ref{sec:Case}, we simplify the definition of the zones by conducting a statistical study on which lines get congested.

\subsection{FlexRequest Price Discovery}  \label{DSO_price_discovery}

In this work, the focus is on the quantity determination rather than on the price determination for the \emph{FlexRequests}. However, we here sketch out one possible approach for determining the DSO's willingness to pay for flexibility.

\emph{FlexOffers} can result from the time- and temperature-variant opportunity cost of FSPs. We want to avoid bias from assumptions on the \emph{FlexOffer} prices which are difficult to predict and vary continuously. Therefore, we focus on determining the actual power reserves (flexibility) needed and the DSO's willingness to pay. 
We assume that the price of \emph{FlexRequests} ${\lambda}_r$ represents an upper bound of the DSO's willingness to pay for flexibility rather than grid upgrades. It can thus be determined ex-ante by the DSO. 
A straightforward method of computing the DSO's willingness to pay ${\lambda}_r$ is to scale the long-term costs of network reinforcements for the DSO with the required flexibility

\begin{equation}
     {\lambda}_r = \frac{ C^\text{Inv} - C^\text{NoInv}}{ P^\text{Flex}} \quad\quad \left[\frac{\$}{\text{kW}}\right]
  \label{eq:problemN}
\end{equation}
where $C^\text{Inv}$ are the total costs of network investments for the DSO if no flexibility was available, and $C^\text{NoInv}$ are the total network investment costs if sufficient inexpensive flexibility $P^\text{Flex}$ was available. $C^\text{Inv},C^\text{NoInv}$ and $P^\text{Flex}$ are the results of two separate investment planning problems that must be solved for the same investment horizon.

\section{Market Clearing} \label{sec:marketclearing}
In this section, we introduce two market-clearing models. The first clears the \emph{FlexRequests} created with the method given in the previous section. The second is a stochastic market clearing for which no \emph{FlexRequest} is submitted as they are implicitly considered in the model. This second market clearing will be used as a benchmark. 

\subsection{Deterministic Market with Explicit FlexRequest} \label{sec:Mdeter}
We assume here that the \emph{FlexRequests}, as determined in Section~\ref{sec:Mcreation} are submitted to the simplest flexibility market architecture, similar to the ancillary services markets currently used on the transmission level. Such a market matches flexibility offers with flexibility requests neglecting all network constraints and not considering uncertainty; the network constraints and the uncertainty have, instead, been already accounted for during the \emph{FlexRequests} creation.

This deterministic market clearing is formulated as
\begin{subequations}
\allowdisplaybreaks
\begin{align}
    \label{eq:DCObj} \underset{p}{\min} & \sum_{t \in \mathcal{T}} \left ( \sum_{o \in \mathcal{O}_t} \lambda_o p_{o} - \sum_{r \in \mathcal{R}_t} \lambda_r p_{r} \right )\\
    \label{eq:DCup} \text{s.t.} & \sum_{o\in \mathcal{O}_{z,t}^+}  p_{o} = \sum_{r\in \mathcal{R}_{z,t}^+}  p_{r}, & \forall z \in \mathcal{Z}, \ \forall t \in \mathcal{T}, \\
    \label{eq:DCdown} & \sum_{o\in \mathcal{O}_{z,t}^-}  p_{o} = \sum_{r\in \mathcal{R}_{z,t}^-}  p_{r}, & \forall z \in \mathcal{Z}, \ \forall t \in \mathcal{T}, \\
    \label{eq:DCboundO} & 0 \le p_{o} \le \overline{P}_{o}, & \forall o\in \mathcal{O}, \\
    \label{eq:DCboundR} & 0 \le p_{r} \le \overline{P}_{r}, & \forall r\in \mathcal{R}.
\end{align}
\end{subequations}
Here, $p$ is the quantity accepted from the corresponding offer or request. The sets $\mathcal{O}_{z,t}$ and $\mathcal{R}_{z,t}$ group respectively the offers and the requests submitted per zone $z$ and time period $t$. The set $\mathcal{Z}$ groups the zones and $\mathcal{T}$ the time periods. The objective function~\eqref{eq:DCObj} is to maximize the social welfare with $\lambda$ the price of the corresponding bid. The constraints~\eqref{eq:DCup}--\eqref{eq:DCdown} ensure that matches are only possible in a given zone, respectively for upwards and downwards bids. Equations~\eqref{eq:DCboundO}--\eqref{eq:DCboundR} give the bounds of the quantity accepted, which is limited by the quantity of the bid $\overline{P}$; sets $\mathcal{O}$ and $\mathcal{R}$ include all offer and request bids for all zones and time periods.

Such a market ensures high transparency for all market participants due to its simplicity, is compatible with the current market regulations as the flexibility market operator does not need access to distribution network data, and we expect that it is the most probable to be implemented in the near future. 

\subsection{Stochastic Market with Implicit FlexRequest}
\label{sec:Mstoch}
We benchmark the performance of the market in \cref{sec:Mdeter} against a network-aware stochastic market that internalizes uncertainty directly in the market clearing. From a theoretical perspective, despite being difficult to exist in practice at the moment, such a market shall determine the most cost-efficient flexibility procurement; that is why it is used as a benchmark. Here again, we use a chance-constrained approach:
\begin{subequations}
\allowdisplaybreaks
\begin{align}
   \nonumber \min_\mathbf{x} & \sum_{t \in \mathcal{T}} \left ( \sum_{o \in \mathcal{O}_t} \lambda_{o}^\text{O} p_{o}^\text{O}+ \mathbb{E} \left [ \Delta t \sum_{n \in \mathcal{N}} \lambda_{n,t}^\text{A} (\tilde{p}_{n,t}^\text{A+} + \tilde{p}_{n,t}^\text{A-}) \right. \right. \\
   \label{eq:SCobj} & \quad \quad \quad \quad \quad \quad \quad \quad \quad \quad \left. \left. +\lambda_{n,t}^\text{NS} \tilde{p}_{n,t}^\text{NS} + \lambda_{n,t}^\text{C} \tilde{p}_{n,t}^\text{C} \right ] \right ) \\
    \nonumber \text{s.t.} & \sum_{jn \in \mathcal{L}} P_{jn,t} - \sum_{nk \in \mathcal{L}} P_{nk,t} = P_{n,t}^\text{inj}+P_{n,t}^\text{U} \\
    \label{eq:SCbalP} &  \quad +p_{n,t}^\text{A} -p_{n,t}^\text{C}+p_{n,t}^\text{NS}, \quad  \quad \quad  \ \forall n \in \mathcal{N},\ \forall t \in \mathcal{T}, \\
    \nonumber & \sum_{jn \in \mathcal{L}} Q_{jn,t} - \sum_{nk \in \mathcal{L}} Q_{nk,t} = K ( P_{n,t}^\text{inj}+P_{n,t}^\text{U} \\
    \label{eq:SCbalQ} &  \quad +p_{n,t}^\text{A} -p_{n,t}^\text{C}+p_{n,t}^\text{NS}), \, \; \quad \quad \ \forall n \in \mathcal{N}, \ \forall t \in \mathcal{T}, \\
    \nonumber &  u_{j,t} = u_{i,t} \\
    \label{eq:SCu} & - 2(R_{ij}P_{ij,t}+X_{ij}Q_{ij,t}),  \, \; \; \quad \  \forall ij \in \mathcal{L}, \ \forall t \in \mathcal{T}, \\
    \label{eq:SCpm} & p_{n,t}^\text{O-} = \sum_{o \in \mathcal{O}_{n,t}^-} p_o^\text{O}, \; \; \, \quad \qquad \qquad \forall n \in \mathcal{N}, \ \forall t \in \mathcal{T}, \\
    \label{eq:SCpm2} & p_{n,t}^\text{O+} = \sum_{o \in \mathcal{O}_{n,t}^+} p_o^\text{O}, \; \; \, \quad \qquad \qquad  \forall n \in \mathcal{N}, \ \forall t \in \mathcal{T}, \\
    \label{eq:SCoff} & 0\leq p_{o}^\text{O} \leq P_{o}^\text{O}, \qquad \qquad \qquad \qquad \qquad \quad  \forall o \in \mathcal{O},\\
    \nonumber  & \mathbb{P}(\tilde{P}_{ij,t}^2 + \tilde{Q}_{ij,t}^2 \leq \overline{S}_{ij}^2) \\ 
    \label{eq:SCcc1} & \quad \geq 1 - \epsilon_\text{S}, \; \; \, \quad \qquad \qquad \qquad \forall ij \in \mathcal{L}, \ \forall t \in \mathcal{T},\\
    \label{eq:SCcc2} & \mathbb{P}(\underline{V}_n^2 \leq \tilde{u}_{n,t}) \geq 1 - \epsilon_\text{V}, \; \; \, \qquad \forall n \in \mathcal{N}, \ \forall t \in \mathcal{T}, \\
    \label{eq:SCcc3} & \mathbb{P}(\tilde{u}_{n,t} \leq \overline{V}_n^2) \geq 1 - \epsilon_\text{V}, \; \; \, \qquad  \forall n \in \mathcal{N}, \ \forall t \in \mathcal{T}, \\
    \label{eq:SCcc4} & \mathbb{P}(- p_{n,t}^\text{O-} \leq \tilde{p}_{n,t}^\text{A}) \geq 1 - \epsilon_A, \;  \; \, \quad  \forall n \in \mathcal{N}, \ \forall t \in \mathcal{T}, \\
    \label{eq:SCcc5} & \mathbb{P}(\tilde{p}_{n,t}^\text{A} \leq p_{n,t}^\text{O+}) \geq 1 - \epsilon_\text{A}, \; \; \qquad  \forall n \in \mathcal{N}, \ \forall t \in \mathcal{T}, \\
    \label{eq:SCcc6} & \mathbb{P}(0 \leq \tilde{p}_{n,t}^\text{C}) \geq 1 - \epsilon_\text{C}, \; \; \, \quad  \qquad   \forall n \in \mathcal{N}, \ \forall t \in \mathcal{T}, \\
    \label{eq:SCcc7} & \mathbb{P}(0 \leq \tilde{p}_{n,t}^\text{NS}) \geq 1 - \epsilon_\text{NS}, \; \quad  \qquad  \forall n \in \mathcal{N}, \ \forall t \in \mathcal{T},
\end{align}
\end{subequations}
where $\mathbf{x}=\{p^\text{O},p^\text{A},p^\text{C},p^\text{NS}, u, P, Q\}$.

This formulation is very similar to the one introduced in Section~\ref{sec:Mcreation}. The main difference is that the offers are included as variables here, in order to be cleared. The objective function~\eqref{eq:SCobj} is to minimize costs for the DSO, as the sum of market costs and expected real-time costs. Here, $\tilde{p}_{o}^\text{O}$ is the quantity accepted from offer $o$ and $\lambda_{o}^\text{O}$ the price of this offer. The real-time costs are the costs associated with activation $\tilde{p}^\text{A}$, load shedding $\tilde{p}^\text{NS}$ and curtailment $\tilde{p}^\text{C}$. The corresponding costs are $\lambda^\text{A}$, $\lambda^\text{NS}$ and $\lambda^\text{C}$. In the objective function, the activation is separated between its positive and negative parts such that $\tilde{p}^\text{A} = \tilde{p}^\text{A+} - \tilde{p}^\text{A-}$, with $\tilde{p}^\text{A+}, \tilde{p}^\text{A-} \geq 0$. This ensures that activation always corresponds to a cost for the DSO.
The multiplication by the duration of the time period $\Delta t$ is needed for coherence between power and energy units. Equations~\eqref{eq:SCbalP}--\eqref{eq:SCu} are defined similarly as in the problem of Section~\ref{sec:Mcreation}.
Equations~\eqref{eq:SCpm}--\eqref{eq:SCpm2} give the total quantity accepted at each node and for each time period, $p_{n,t}^\text{O+}$ in the up direction, and $p_{n,t}^\text{O-}$ in the down direction. The quantity accepted per offer is bounded by the quantity offered $P^\text{O}$ in Equation~\eqref{eq:SCoff}. Equations~\eqref{eq:SCcc1}--\eqref{eq:SCcc7} are chance constraints. Equations~\eqref{eq:SCcc4}--\eqref{eq:SCcc5} correspond to the bounds on the activation of flexibility, which are given by how much flexibility is procured. Equations~\eqref{eq:SCcc6}--\eqref{eq:SCcc7} are the lower bounds on curtailment and load shedding. In those constraints, $\epsilon_\text{S}$, $\epsilon_\text{V}$, $\epsilon_\text{A}$, $\epsilon_\text{C}$ and $\epsilon_\text{NS}$ are used to set the confidence level.

The reformulation of the chance constraints is very similar to the one shown in Section~\ref{sec:Mcreation} and will not be detailed here. In the case of the stochastic market clearing, each of the variables $\tilde{p}^\text{A}$, $\tilde{p}^\text{NS}$ and $\tilde{p}^\text{C}$ is expressed as a linear function of the total forecast error, using coefficients $\alpha^\text{A}$, $\alpha^\text{NS}$ and $\alpha^\text{C}$, such that the nodal $\alpha$ is $\alpha_{n,t}=\sum_{n \in \mathcal{N}} \alpha_{n,t}^\text{A} + \alpha_{n,t}^\text{NS} + \alpha_{n,t}^\text{C} = 0, \ \forall t \in \mathcal{T}$. Balance responsibility can be ensured by constraints similar to those in \ref{sum_bal2}.
The reformulation of the expected costs in the objective function is straightforward with the assumption that the forecast error follows a Gaussian distribution with zero mean. In such a market clearing, the FMO must obtain full knowledge of the underlying network, the setpoints of all resources, and information about the uncertainty (covariance matrix).

\subsection{Role of the DSO \& Computational Tractability} 
In the benchmark, the FMO runs a complex stochastic market clearing with network constraints of the DSO's distribution network data. This requires computational efforts from the FMO, which are commonly covered by market participation fees.

However, if the DSO insists on distribution network data privacy, then the only possible way to clear the distribution-level flexibility market in a network-aware fashion is the approach we are proposing. This also shifts the computational burden away from the FMO and onto the DSO, which is also the entity that has the need for the flexibility service. In that way, the DSO solves the problem in \cref{sec:Mcreation} to create network-aware \emph{FlexRequests} that can easily be handled in a deterministic zonal market by the FMO.

\section{Suboptimality Gap} \label{sec:subopt}
\subsection{Suboptimality Gap of Network-Aware FlexRequest Compared to Stochastic Clearing}
By intuition, the use of network-aware \emph{FlexRequests}, which are requested at specific nodes or zones, pre-assuming the existence of possible congestions, may introduce suboptimality compared to the stochastic-clearing benchmark, which makes no prior assumptions on possible line constraint or voltage violations. This is a common well-identified issue in the transmission-level literature on zonal vs. nodal markets, e.g. \cite{Wang2015}. Nodal markets result in higher social welfare because they include no prior assumption about the designation of the congestion-free zones. Therefore, in this section, we quantify and derive an analytical expression to determine this suboptimality gap. 
In the following, we define ``sufficient liquidity'' as sufficient offer volume $\sum_n p_n^\text{O}$ in the entire network. We define ``sufficient liquidity per node'' as sufficient offer volume $p_n^\text{O}$ at each node, which is offered at the same or lower prices than that of other nodes.
We can distinguish four levels of increasing complexity:
\begin{enumerate}
  \item \emph{No active network constraints and sufficient liquidity}: Here, the market-clearing outcome of the network-aware \emph{FlexRequests} and of the stochastic-clearing benchmark are the same.
  \item \emph{Active network constraints and sufficient liquidity per node}: Here, the market-clearing outcome of both cases are the same.
  \item \emph{Active network constraints, sufficient liquidity in the network, but insufficient liquidity per node}: This is the relevant case; the network-aware \emph{FlexRequests} may result in a suboptimal outcome compared to the stochastic-clearing benchmark.
  \item \emph{Active network constraints and insufficient liquidity in the network}: The market-clearing outcome results in load shedding and/or production curtailment in both cases which is an undesired outcome for local flexibility markets.
\end{enumerate}
In the following, we will evaluate the market-clearing outcomes of level 1 and level 3 in terms of social welfare.

\subsection{Upper Bound on Suboptimality Gap}
The market-clearing problem without active network constraints and sufficient liquidity (level 1) can be described by the unconstrained market-clearing outcome $\mathcal{L}_U$:
\begin{equation}
  \begingroup
  \allowdisplaybreaks
  \mathcal{L}_U=
   \begin{cases}
  	& \underset{p_n^\text{O}}{\max} \; \lambda^\text{R} \sum_n p_n^\text{R} - \sum_{n} \lambda_{n}^\text{O} p_{n}^\text{O} \\
	& \text{s.t.} \;\sum_{n} p_{n}^\text{O} = \sum_n p_n^\text{R}  \;\; \\
    & \;\;\;\;\;\;\lambda_{n}^\text{O} \le \hat{\lambda}^\text{R} \;\;\forall n
    \end{cases}
  \endgroup .
  \label{eq:problemUnconstr}
\end{equation}
The theoretical worst-case market-clearing outcome occurs when there is insufficient liquidity (level 4), and only becomes feasible through load shedding and/or curtailment. This is the worst case since market transactions are constrained by both insufficient \emph{FlexOffers} and active network constraints. Additionally, social welfare is reduced by payments to loads due to load shedding and/or to generators due to curtailment.

The more interesting comparison is with the worst-case feasible market-clearing outcome when there is sufficient liquidity in the whole network, but insufficient liquidity per node (level 3). Thus, the flexibility may not be provided at the same node, and the additional loading may result in congestion or voltage limit violations. These active network constraints may limit the volume of the accepted offer which can be characterized by a permitted share of the offer $a_n p_n^\text{O}$ and a non-permitted share $(1-a_n) p_n^\text{O}$ where $a_n\in(0,1), \, \forall n$.

The market-clearing outcome with stochastic market clearing is given by \cref{eq:problemSC}:
\begin{equation}
  \begingroup
  \allowdisplaybreaks
  \mathcal{L}_{SC}=
   \begin{cases}
  	& \underset{p_n^\text{O}}{\max} \;\lambda^\text{R} \sum_n p_{n}^\text{R} - \sum_{n} a_n \lambda_{n}^\text{O} p_{n}^\text{O} \\
	& \text{s.t.} \;\sum_{n} a_n p_{n}^\text{O} = \sum_n p_n^\text{R}  \;\;\ \\
    & \;\;\;\;\;\;\lambda_{n}^\text{O} \le \hat{\lambda}^\text{R}  \;\;\forall n
    \end{cases}
  \endgroup .
  \label{eq:problemSC}
\end{equation}
Ignoring losses, the upper bound on the suboptimality gap compared to the unconstrained clearing is then given by \cref{eq:BoundSC}:
\begin{equation}
  \Xi_{SC} = 1- \frac{\mathcal{L}_{SC}}{\mathcal{L}_U}.
  \label{eq:BoundSC}
\end{equation}

Similarly, the location-based market clearing outcome with \emph{FlexRequests} is given by \cref{eq:problemFR}:
\begin{equation}
  \begingroup
  \allowdisplaybreaks
  \mathcal{L}_{FR}=
   \begin{cases}
  	& \underset{p_n^\text{O}}{\max} \; \lambda^\text{R} \sum_{n}  p_{n}^\text{R} - \sum_{n} a_n \lambda_{n}^\text{O} p_{n}^\text{O} \\
	& \text{s.t.} \;a_n p_{n}^\text{O} = p_n^\text{R}  \;\;\forall n \\
    & \;\;\;\;\;\;\lambda_{n}^\text{O} \le \hat{\lambda}^\text{R}  \;\;\forall n
    \end{cases}
  \endgroup ,
  \label{eq:problemFR}
\end{equation}
and the upper bound on its suboptimality gap is given by \cref{eq:BoundFRU}:
\begin{equation}
  \Xi_{FR} = 1- \frac{\mathcal{L}_{FR}}{\mathcal{L}_{U}}.
  \label{eq:BoundFRU}
\end{equation}

Furthermore, and most relevant, the upper bound on the suboptimality gap of the use of \emph{FlexRequests} compared to stochastic clearing can be expressed by \cref{eq:BoundFRSC}.
\begin{equation}
  \Xi_{FS} = 1- \frac{\mathcal{L}_{FR}}{\mathcal{L}_{SC}}.
  \label{eq:BoundFRSC}
\end{equation}
In Section \ref{sec:Case}, we show some results regarding suboptimality gaps of the different methods.

\section{Case Study}\label{sec:Case}
This section presents two case studies, one on a 15-bus system with publicly available data, and a second on a segment of a real German distribution network.

\subsection{15-bus Test Case}
We evaluate our method for the creation of \emph{FlexRequests} on the radial 15-bus system introduced in~\cite{das1995}. Note that the true complexity does not stem from the size of the power system as much as it does from the size of the uncertainty. We consider uncertainty from renewable production, while the rest of the power injections are assumed to be certain for the purpose of this example. Two wind farms are placed in the system and the covariance matrix is evaluated from 1,000 wind forecast error scenarios with zero mean. These are taken from real measurements recorded in Denmark~\cite{pinson2013}. Without loss of generality, we consider that the reserve market is cleared for one time period of one hour, and a power factor $\cos  \phi = 0.95$. Solving for multiple time periods would only require the solution of the same problem an equal number of times, under the assumption that uncertainty at time $t$ and $t+1$ are statistically independent. The exact data used is available online, along with the implemented models~\cite{GitHub}.

Regarding the different prices, we make the following assumptions:
\begin{itemize}
    \item Curtailment costs are set at 60€/MWh.
    \item Load shedding costs at 200€/MWh.
    \item Activation costs are assumed to be equal to 0 since we consider that the flexibility providers have no operational costs.
    \item \emph{FlexRequests} have a price of 70€/MW for up- and 40€/MW for down-regulation.
    \item Offers are priced randomly between 25€/MW and 35€/MW.
\end{itemize}

For all the different chance constraints, a violation of 5\% is permitted.

The following are solved:
\begin{itemize}
    \item \emph{FlexRequest} creation with the chance-constrained model.
    \item Deterministic market clearing to match the created \emph{FlexRequests} with \emph{FlexOffers} submitted in the market.
    \item Stochastic market clearing as a benchmark to compare with the deterministic market clearing which considered the created \emph{FlexRequest}.
    \item Out-of-sample Monte Carlo analysis to evaluate the violation of chance constraints for both models. 
    \item Real-time dispatch (with \emph{LinDistFlow}) to evaluate the real-time costs after activation.
\end{itemize}

The models are implemented in Python, with the Pyomo Kernel library \cite{pyomo} and solved with Mosek for the chance-constrained problems and with the Gurobi library \cite{gurobi} for the deterministic market clearing.

\subsection{Results for the 15-bus System}

The performance of the model for the creation of the \emph{FlexRequests} is evaluated over two axes. The first is the definition of clearing zones by the DSO. The most extreme case is to have a clearing per node, which corresponds to either not having the possibility to define zones or to the DSO considering that all lines are likely to get congested such that any deviation shall be balanced through flexibility procurement at the same node it is created. This case will be identified as ``Per-Node". For the second case, we conducted a statistical study to identify which lines can get congested and defined zones accordingly. This will be further referred to as ``4-Zones". Finally, we define a test case that is between those two, with more zones, and it will be called ``8-Zones".

The second parameter is the overall liquidity of the market. The offers are designed following three levels of liquidity: high, medium, and low. For the high-liquidity scenario, all buses except the reference bus are offering a high level of flexibility in both directions. This is then reduced to obtain the medium and low scenarios.

\begin{figure}[t]
  \centering
    {\includegraphics[width=0.8\linewidth]{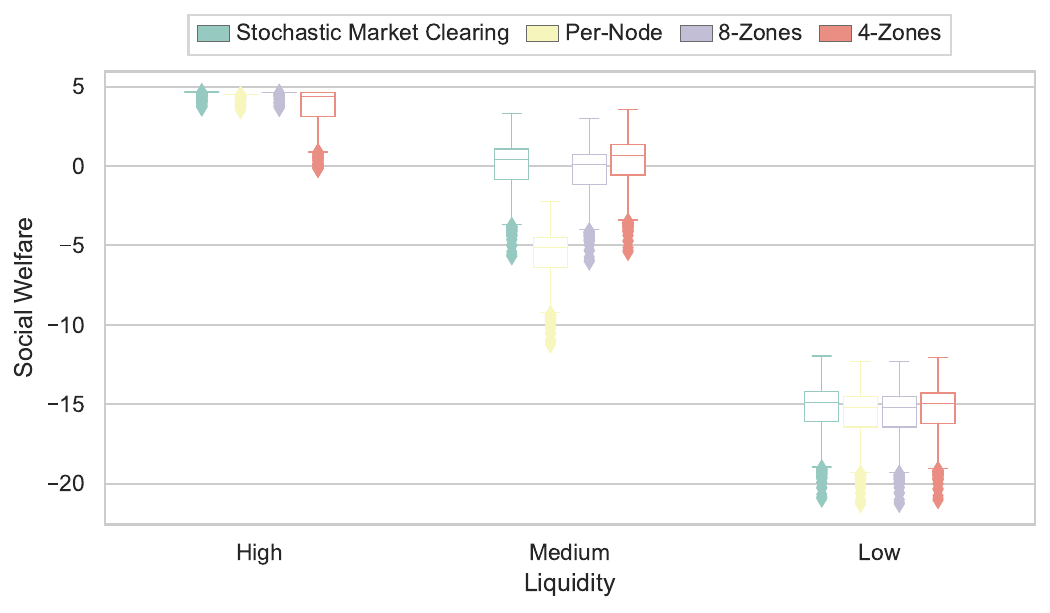}}
  \vspace{-0.5cm}
  \caption{Social welfare (€) for three levels of liquidity and comparing stochastic market clearing to deterministic market clearing with different ways to define clearing zones, for the 15-bus test case}
  \label{fig:SW}
\end{figure}

\begin{figure}[tb]
    \centering
    \includegraphics[width=0.85\linewidth]{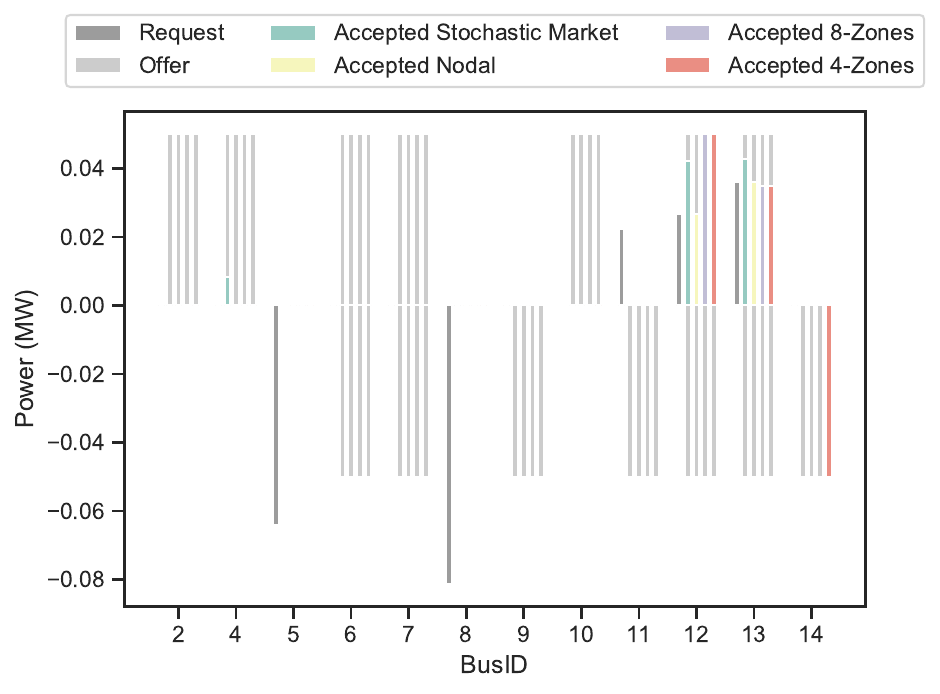}
    \caption{Submitted offers and requests per bus in the case of medium liquidity. The upper part of the graph corresponds to bids up and the lower part to bids down. The filling of the offers bar shows the quantity accepted after the corresponding market clearing.}
    \label{fig:bids}
\end{figure}

After running those, we compare the resulting social welfare. It is calculated as
\begin{equation}
\small
    \text{SW} = (\lambda^\text{R+}-\mathbf{\lambda}^\text{O+}) \mathbf{p}^\text{f+} + (\lambda^\text{R-}-\mathbf{\lambda}^\text{O-}) \mathbf{p}^\text{f-} - \lambda^\text{NS} \mathbf{p}^\text{NS} - \lambda^\text{C} \mathbf{p}^\text{C}.
    \label{eq:SW}
\end{equation}
The first part of the social welfare is linked to procurement and is calculated once. It is the difference between the price of accepted requests and offers multiplied by the procured quantity. 
The second part corresponds to the costs resulting from real-time dispatch, which are due to load shedding and curtailment. Those are calculated over 2,000 scenarios of wind realization, which are different from the ones used to evaluate the covariance matrix. 
The social welfare for the different test cases is shown in Figure~\ref{fig:SW}. Note that the values can be negative, which is due to the fact that there is some load shedding, for which the costs are considerably higher than the other prices.
The more offers are available, the higher the social welfare. This is due to both the social welfare gained from the matches and the reduction of real-time redispatching costs. The deterministic market with no zones (i.e. Per-Node) can perform significantly worse than the others, as we see here in the medium liquidity. It is consistent with the fact that it restricts the possibility of matches. Since less flexibility is procured, the real-time costs, linked to load shedding in particular, are also higher. The performance of the deterministic market with zones is close to the one of the stochastic market in terms of social welfare. With fewer, larger clearing zones (4-Zones), more matches are allowed, so with similar redispatch costs, this can give an advantage to the deterministic market in terms of social welfare, but it means that some flexibility is not used. This happens when the zones are defined too large, which means that there exists congestions inside of the zone that prevent the activation of the procured flexibility. In Figure~\ref{fig:bids}, we can see that with 8-Zones, some downward flexibility is procured in bus 14, while the stochastic market does not have matches in this direction because the location of the offers is not helpful. On this plot, we can also clearly see that with the zonal market, the procured flexibility corresponds to the quantity requested in buses 12 and 13, while the stochastic and zonal markets procure more at these buses.

\begin{figure*}[t]
  \centering
    {\includegraphics[width=0.70\linewidth]{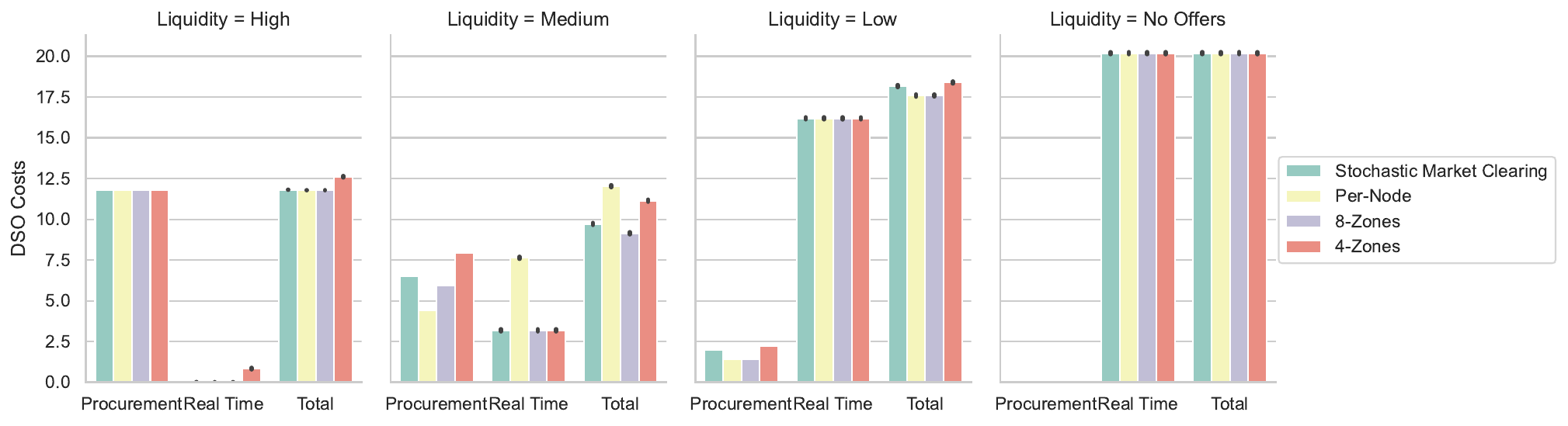}}
  \vspace{-0.5cm}
  \caption{DSO costs (€) for three levels of liquidity and comparing stochastic market clearing to deterministic market clearing with different ways to define clearing zones. The last graph illustrates the real-time costs without a flexibility market, for the 15-bus test case}
  \label{fig:DSOcosts}
\end{figure*}

\begin{figure*}[t]
  \centering
    {\includegraphics[width=0.70\linewidth]{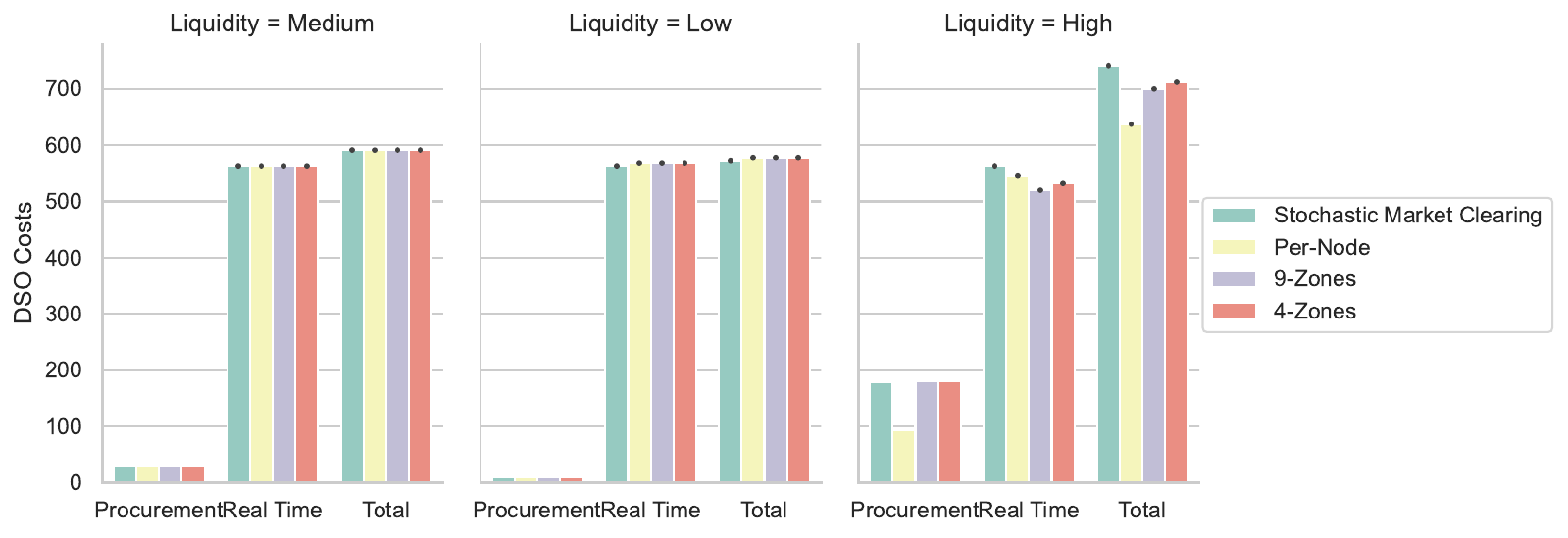}}
  \vspace{-0.5cm}
  \caption{DSO costs (€) for three levels of liquidity and comparing stochastic market clearing to deterministic market clearing with different ways to define clearing zones, for the bnNETZE 81-bus test case}
  \label{fig:DSOcosts_81}
\end{figure*}

In Figure~\ref{fig:DSOcosts}, the DSO costs in the different situations are compared. The procurement costs for the DSO depend on the payment scheme chosen in each market. To keep the market-clearing models general in this regard, we consider that the DSO have to pay what they bid, which would be a higher bound on the price paid. 
The graph on the right (Liquidity=No Offers) illustrates what the costs would be without a flexibility market. The costs always decrease with the flexibility market, compared to not having one (in the latest plot).
With the clearing per node, the quantity matched at each node can be reduced, which gives lower procurement costs for medium and low liquidity here. In the case of medium liquidity, the quantity procured at a critical bus is much lower compared with the zonal market, where flexibility requested at a critical bus can be reserved at several surrounding nodes. As a consequence, real-time costs are much higher in that case. It is worth noting here that since the \emph{FlexRequest} creation model is based on an optimization algorithm, the optimization solver by default returns only a single possible combination for the location of the requests; there could be, however, more than one solution. Having zones is one way to deal with this limitation, as we can see in this example. The level of liquidity plays a role in the allocation of the costs between procurement and real time.

The performance of the chance-constrained models was also analyzed and the maximum probability of violation over all the chance constraints was evaluated in an out-of-sample analysis over 2,000 scenarios. It was found to be equal to 4\% for both models, which is in line with the limit of 5\% that we have defined for all chance constraints.

In terms of runtime, Table \ref{tab:results_runtime} shows that the deterministic market clearing is much faster than the stochastic market clearing, due to the fact that the computational burden is moved to the \emph{FlexRequest} creation in the second situation. Table \ref{tab:results_so_gap} gives results in terms of suboptimality. Those are obtained by placing ourselves in the situation described in Section \ref{sec:subopt}. We compare the case with no congestion to stochastic market clearing and to nodal deterministic market clearing, for the high-liquidity scenario. The gaps in the case with no congestion are due to the fact that there are congestions. More importantly, we see that the gap between the nodal deterministic market clearing and the stochastic market clearing is only 3.2\% in this case.

\begin{table}[bt]
\caption{Results in terms of runtime (in average over all case studies)}
\centering
\begin{tabular}{l|r|r|}
\cline{2-3}
\textbf{} & \multicolumn{1}{c|}{\textbf{15-bus}} & \multicolumn{1}{c|}{\textbf{bnNETZE}} \\ \hline
\multicolumn{1}{|l|}{Stochastic}           & 2.75 s & 1.33 h \\
\multicolumn{1}{|l|}{Deterministic}        & 0.07 s & 1.04 s \\
\multicolumn{1}{|l|}{FlexRequest creation} & 2.88 s & 1.34 h \\ \hline
\end{tabular}
\label{tab:results_runtime}
\end{table}

\begin{table}[bt]
\caption{Suboptimality gap for both test systems}
\centering
\begin{tabular}{|ll|r|r|}
\cline{3-4}
\multicolumn{1}{c}{} &
  \textbf{} &
  \multicolumn{1}{c|}{\textbf{15-bus}} &
  \multicolumn{1}{c|}{\textbf{bnNETZE}} \\ \hline
 \multirow{3}{*}{\textbf{Market clearing outcome}}  & $\mathcal{L}_{U} $ \eqref{eq:problemUnconstr} & 5.83 € & 103.04 €\\
    & $\mathcal{L}_{SC} $ \eqref{eq:problemSC} &  4.66 €  & 90.83 €  \\
    & $\mathcal{L}_{FR} $ \eqref{eq:problemFR} &  4.51 €  & 99.00 € \\ \hline
\multirow{3}{*}{\textbf{Suboptimality gap}} & $\Xi_{SC}$ \eqref{eq:BoundSC} & 20.1 \% & 3.92 \% \\
              &  $\Xi_{FR}$ \eqref{eq:BoundFRU} & 22.6 \% & 11.9 \% \\
              & $\Xi_{FS}$ \eqref{eq:BoundFRSC} & 3.2 \% & 8.25 \% \\ \hline
\end{tabular}
\label{tab:results_so_gap}
\end{table}

\subsection{Analysis on a Real Distribution Grid}
\begin{figure}[t]
  \centering
    {\includegraphics[width=1\linewidth]{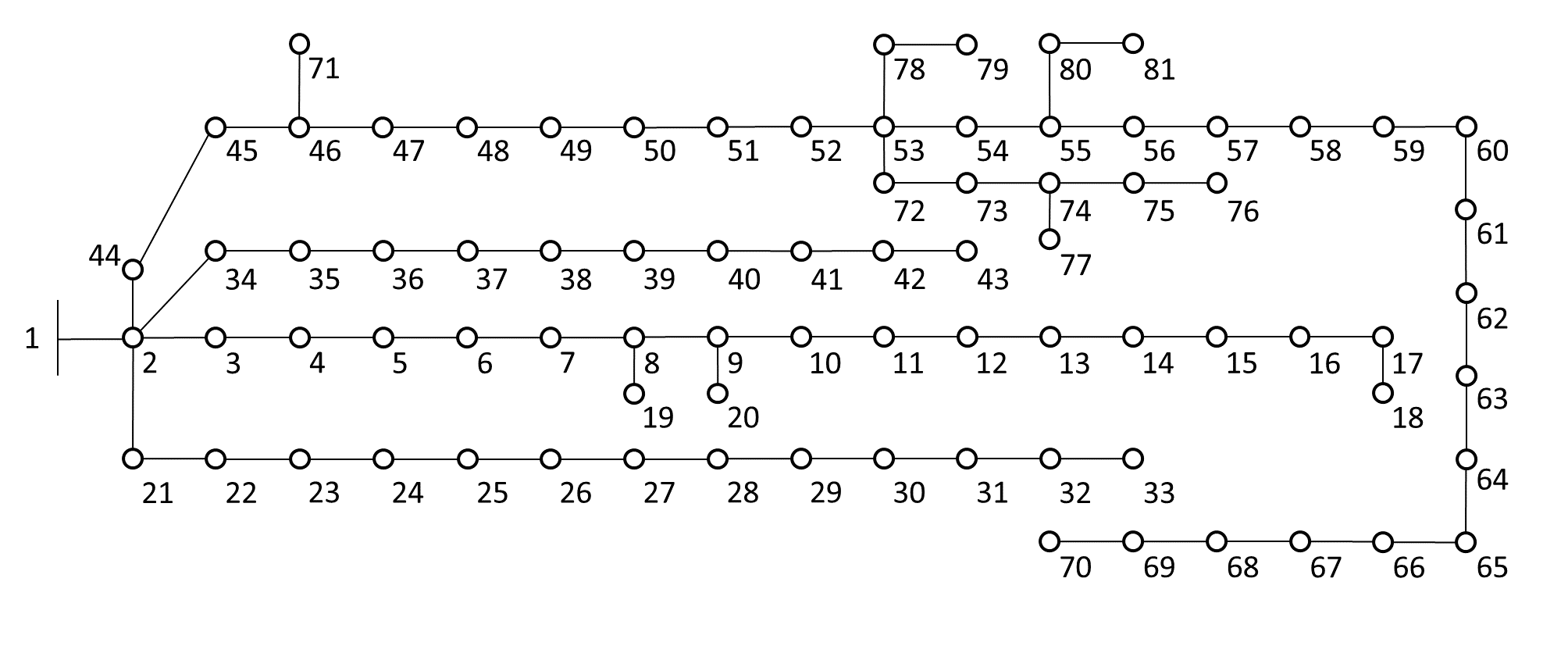}}
  \vspace{-0.5cm}
  \caption{Representation of the bnNETZE 81-bus network}
  \label{fig:network}
\end{figure}

\begin{figure}[t]
  \centering
    {\includegraphics[width=0.8\linewidth]{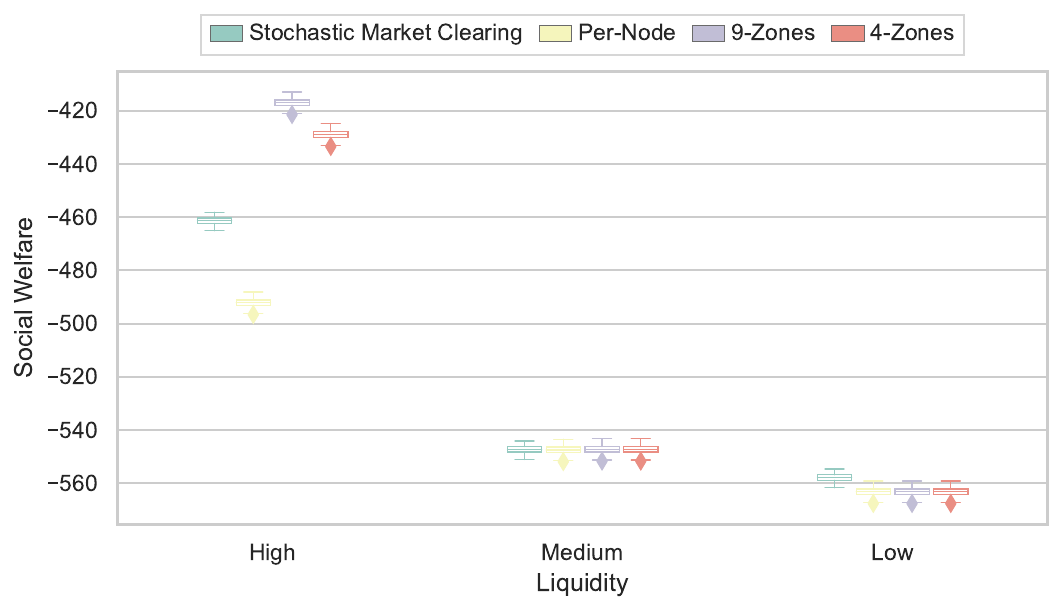}}
  \vspace{-0.5cm}
  \caption{Social welfare (€) for three levels of liquidity and comparing stochastic market clearing to deterministic market clearing with different ways to define clearing zones, for the bnNETZE 81-bus test case}
  \label{fig:SW_81}
\end{figure}

In order to show scalability, we apply the same methods to a radial 81-bus German distribution network, which was provided by the DSO bnNETZE. It is presented in Figure~\ref{fig:network}. The 81 bus 0.4kV distribution system is connected to the substation through a 20kV line and has a total maximum loading of 3.5MW, and three wind farms with a total rated capacity of 500kW. The other parameters are defined similarly as in the previous study, and in particular, the same data is used to represent the uncertainty. Similarly, we define different types of zones for the deterministic market clearing, identified as ``4-Zones" and ``9-Zones". 
The results in terms of social welfare and costs for the DSO are shown in Figure~\ref{fig:SW_81} and~\ref{fig:DSOcosts_81} respectively. The conclusions of the analysis for the 15-bus system are still valid here. The particularity of this case is that the system is heavily loaded. The real-time costs are very high in all the scenarios tested, due to load shedding. The stochastic market clearing is performing better, as expected, but the zonal deterministic markets are comparatively close.

As shown in Table \ref{tab:results_runtime}, the running time difference between deterministic market clearing and stochastic market clearing is even more pronounced in this case, as the deterministic market clearing is more than 4,600 times faster. Moreover, the suboptimality gap between the \emph{nodal} deterministic market clearing and the stochastic market clearing is still reasonable, with 8.25\% in this case.

\subsection{Discussion}
This section summarizes the advantages and disadvantages of the proposed privacy-preserving flexibility market compared to the stochastic market clearing. It is based on the points presented in Section I and the results from the two case studies. In terms of preserving privacy, the proposed model inherently has this property compared to the stochastic market model. Privacy preservation comes at the expense of reduced social welfare. We provide upper bounds on the suboptimality gap of the deterministic flexibility market compared to the stochastic one, ranging from 3.2 \% for the small network to 8.25 \% for the real-size network. Gaps of this magnitude may still be acceptable, especially when considering the benefits of market transparency and computational tractability. Moreover, the case studies show that this gap can be reduced by a proper definition of zones for the market clearing. We also observed that the results for social welfare are greatly influenced by the market liquidity, i.e., the number of \emph{FlexOffers} at each location, which is an external parameter. Similarly, the costs for the DSO are impacted by the definition of zones and by market liquidity. We show that social welfare and DSO costs are actually comparable for the stochastic and the deterministic network-aware market when zones are defined. Another considerable advantage of the proposed flexibility market model is that it is cleared in the range of seconds for both cases, whereas the stochastic market requires at least an hour for the realistic test case.

\section{Conclusion}\label{sec:Conc}
Although there is a lot of discussion and research literature on the design of efficient distribution-level flexibility markets, little attention has been so far paid to the design of the \emph{FlexRequest} itself. In this paper, we provided a framework for DSOs to create a network-aware \emph{FlexRequest} as a tool that ensures compatibility with current legislation. We compared its efficacy against a stochastic benchmark, showing results for a real-world distribution network. We showed how the \emph{FlexRequests} can be generated with a chance-constrained approach and cleared in a very simple deterministic market, which includes the possibility to define zones. The results of flexibility procurement and real-time dispatch were compared for different levels of liquidity on the offer side. We saw that the definition of bidding zones by the DSO can become critical and that by enforcing the procurement of flexibility at the same node it is requested, the costs for the DSO could increase significantly. With properly defined zones, the results obtained in terms of social welfare and costs for the DSO appeared to be comparable to those obtained with stochastic market clearing. 

In future work, we will focus on the determination of the price of the requests and the definition of the bidding zones for the DSO. The model could be extended to have \emph{FlexRequests} which also consider contingencies. Moreover, the role of reactive power could be detailed and expanded, in particular to consider reactive power bids. Finally, we made the assumption of a radial network, and it could be worth generalizing to any network architecture.

\setcounter{table}{0}
\renewcommand{\thetable}{A\arabic{table}}
\setcounter{figure}{0}
\renewcommand{\thefigure}{A\arabic{figure}}

\bibliographystyle{IEEEtran}
\bibliography{Bib}

\end{document}